\documentclass[11pt]{article}
\usepackage{graphicx, colordvi}
\usepackage{amsfonts}
\usepackage{amssymb}
\usepackage{amsthm,cite,color}
\usepackage{dsfont}
\usepackage{epsfig}
\usepackage{mathrsfs}
\usepackage{amsfonts}
\usepackage{amssymb}
\usepackage{amsmath}
\usepackage{amssymb,amsfonts,amsmath,amsthm,cite,color}
\usepackage{dsfont}
\usepackage{epsfig}
\usepackage{mathrsfs}
\usepackage{longtable}
\usepackage{hyperref}
\hypersetup{
colorlinks=true,
linkcolor=cyan,
filecolor=blue,
urlcolor=red,
citecolor=green,
}

\newtheorem{theo}{Theorem}
\newtheorem{rem}{Remark}
\newtheorem{prop}[theo]{Proposition}

\newtheorem{defi}[theo]{Definition}
\newtheorem{lem}[theo]{Lemma}

\makeatletter \@addtoreset{equation}{section}
\@addtoreset{theo}{section} \makeatother

\newcommand{\bN} { {\mathbb{N}}}
\newcommand{\bC} { {\mathbb{C}}}

\newcommand{\bZ} { {\mathbb{Z}}}

\DeclareMathOperator{\lc}{lc}
\newcommand{\la} { {\langle}}
\newcommand{\ra} { {\rangle}}

\def\qed{\hfill \rule{4pt}{7pt}}
\def\pf{\noindent {\it Proof.} }

\begin{document}
\begin{center}

 {\large \bf $q$-Rational Reduction and $q$-Analogues of Series for $\pi$}
\end{center}
%
%
\begin{center}
{  Rong-Hua Wang}$^{1}$ and {Michael X.X. Zhong}$^{2}$

   $^1$School of Mathematical Sciences\\
   Tiangong University \\
   Tianjin 300387, P.R. China \\
   wangronghua@tiangong.edu.cn \\[10pt]

   $^2$School of Science\\
   Tianjin University of Technology \\
   Tianjin 300384, P.R. China\\
   zhong.m@tjut.edu.cn

\end{center}

\vskip 6mm \noindent {\bf Abstract.}

In this paper, we present a $q$-analogue of the polynomial reduction which was originally developed for hypergeometric terms.
Using the $q$-Gosper representation, we describe the structure of rational functions that are summable when multiplied with a given $q$-hypergeometric term.
The structure theorem enables us to generalize the $q$-polynomial reduction to the rational case, which can be used in the automatic proof and discovery of $q$-identities.
As applications, several $q$-analogues of series for $\pi$ are presented.

\noindent {\bf Keywords}: summability; rational reduction; series for $\pi$.

\section{Introduction}

In 1971, Abramov \cite{Abramov1971} described an algorithm to decide whether an indefinite rational summation $\sum_{k=0}^{n}t_k$ has a closed form.
The problem is equivalent to determining whether the rational function $t_k$ can be written as the difference of another rational function, that is, the \emph{summability} of $t_k$.
Later, Abramov \cite{Abramov1975} provided a method of separating the rational component of the solution of a first-order linear recurrence relation with a rational right side.
In 1978, Gosper \cite{Gosper1978} showed that the summability problem of hypergeometric terms can be reduced to finding rational solutions of first-order difference equations with polynomial coefficients.
Gosper's algorithm is a powerful tool in dealing with different kinds of indefinite summations.

Based on Gosper's algorithm, Zeilberger \cite{Zeilberger1990c, Zeilberger1991} gave a fast algorithm for proving terminating hypergeometric identities and also formulated a vision for doing creative telescoping to holonomic functions \cite{Zeil1990}.
The $q$-analogues of Gosper's algorithm and Zeilberger's algorithm had been developed by Koornwinder \cite{Koornwinder1993}.
Those algorithms have been widely used in combinatorics, number theory and mathematical physics.
At the same time, the mechanical proof of combinatorial identities is still a subject of ongoing research, one can consult \cite{ChenKauers2017} for a selection of open problems in this context.

Recall that a hypergeometric term $t_k$ is a function of $k$ such that the shift quotient $t_{k+1}/t_k$ is a rational function of $k$.
Apparently all nonzero rational functions are hypergeometric.
The summability of rational functions $t_k$ was studied by Abramov \cite{Abramov1971, Abramov1995b}.
The key idea is to decompose a given rational function $t_k$ as
$t_k=\Delta_k (s_k)+r_k$, where $\Delta_k$ is the difference operator with respect to $k$, $s_k$ and $r_k$ are rational functions with the denominator of $r_k$ being shift-free.
Such a decomposition is called a minimal decomposition for $t_k$, and $t_k$ is summable if and only if $r_k$ is zero.
Minimal decompositions of hypergeometric terms and $q$-hypergeometric terms were presented by Abramov and Petkov{\v{s}}ek \cite{AP2001} and Du et al. \cite{DHL2018} respectively.

In general, we care about the decomposition of $t_k$ for two reasons.
First, the minimal decomposition not only tells us whether $t_k$ is summable but also what is the essentially non-summable part when $t_k$ is not summable.
Second, once a decomposition is given, identities or congruence properties for $t_k$ may be deduced.

In 2015, Chen et al. \cite{CHKL2015} presented a modified Abramov--Petkov{\v{s}}ek reduction, which is more efficient than the original reduction and can be used to compute minimal telescopers for bivariate hypergeometric terms.
For any given hypergeometric term $T$,
they first do the multiplicative decomposition $T=SH$, where $S$ is a rational function and $H$ is another hypergeometric term whose shift quotient is shift reduced.
Let $K$ be the shift quotient of $H$.
A crucial step of the modified Abramov--Petkov{\v{s}}ek reduction was the introduction of polynomial reduction with respect to $K$.

Recently, Hou, Mu and Zeilberger~\cite{HouMuZeil2021} introduced a polynomial reduction process for any hypergeometric term (no multiplicative decomposition was needed).
More precisely, they focused on terms of the form $f(k)t_k$, where $f(k)$ is a polynomial and $t_k$ is any hypergeometric term.
Their basic idea is to decompose $f(k)t_k$ as
$
f(k)t_k=\Delta_k(g(k)t_k)+h(k)t_k,
$
where $g(k),h(k)$ are polynomials and the degree of $h(k)$ is bounded.
As applications, they derived two infinite families of supercongruences.
Later on, Hou and Li \cite{HouLi2021} considered the Gosper summability of $r(k)t_k$, where $r(k)$ is a rational function and $t_k$ is hypergeometric, and provided an upper bound and a lower bound on the degree of the numerator of $r(k)$.

We notice that the above polynomial reduction, given by Hou, Mu and Zeilberger, can be generalized to the $q$-case.
Combined with the $q$-Gosper algorithm, we will also characterise when the product of a rational function and a $q$-hypergeometric term can be summable.
Motivated by the corresponding structure theorem, we present a $q$-rational reduction, which can be applied to generate infinitely many new identities.
Our $q$-rational reduction is free of multiplicative decomposition, which differs from the one presented by Du et al. in \cite{DHL2018}.
Moreover, the denominator of the rational part can be assigned in advance.

The paper is organized as follows.
In Section \ref{sec:q poly red}, we introduce the so-called $q$-polynomial reduction.
Section \ref{sec:q rat red} is devoted to determining the structure of $r(x)$ when $r(q^k)t_k$ is summable, which leads to the discovery of $q$-rational reduction.
How to generate new identities, especially those $q$-analogues of the Ramanujan-type series for $\pi$, is presented in Section \ref{sec: application}.

\section{$q$-Polynomial Reduction}\label{sec:q poly red}
Let $\bC$ be the field of complex numbers and $\bC[x]$ the ring of polynomials in $x$ over $\bC$.
Suppose $q\in\bC\setminus\{0\}$ and $|q|<1$.
Then an expression $t_k$ is called a \emph{$q$-hypergeometric} term if the shift quotient is a rational function in $q^k$,
that is, there exist $a(x),b(x)\in\bC[x]$ such that
\begin{equation}\label{eq:q-hyper}
\frac{t_{k+1}}{t_k}=\frac{a(q^{ k})}{b(q^{k})}.
\end{equation}
For any polynomial $g(x)\in\bC[x]$, it is straightforward to check that
\begin{equation}\label{eq:difference equation}
\left(a(q^k)g(q^{k+1})-b(q^{k-1})g(q^{k})\right)t_k
=\Delta_k(b(q^{k-1})g(q^{k})t_k),
\end{equation}
where $\Delta_k$ is the difference operator with respect to $k$.

Given a pair of polynomials $a(x),b(x)\in \bC[x]$, define
\begin{equation}\label{eq:Sab}
S_{a,b}:=\{a(x)g(qx)-b(q^{-1}x)g(x)\mid g(x)\in\bC[x]\}.
\end{equation}
It is clear that $S_{a,b}$ is a subspace of $\bC[x]$, called the $q$-difference space corresponding to $(a(x),b(x))$.

Next we will characterise the degrees of the polynomials in $S_{a,b}$.
Let $\bZ,\bN$ denote the sets of integers and nonnegative integers respectively.
For fixed $a(x),b(x)\in\bC[x]$ with $a(x)b(x)\neq 0$, denote
\begin{equation}\label{eq:d}
 d:=\max\{\deg a(x),\deg b(x)\},
\end{equation}
where $\deg p(x)$ denotes the degree of the polynomial $p(x)$ in $x$.

\begin{defi}
Given $a(x),b(x)\in\bC[x]$ with $d$ defined by \eqref{eq:d}, if
\begin{equation}\label{eq:m_0}
\deg a(x)=\deg b(x) \text{ and } \frac{\lc b(q^{-1}x)}{\lc a(x)}=q^{m_0} \text{ for some } m_0\in\bN,
\end{equation}
where $\lc p(x)$ denotes the leading coefficient of $p(x)$ in $x$,
then the pair $(a(x),b(x))$ is called $\emph{degenerated}$.
\end{defi}
When $\frac{t_{k+1}}{t_k}=\frac{a(q^{k})}{b(q^{k})}$ and $(a(x),b(x))$ is degenerated, we may also say $t_k$ is degenerated when there is no confusion.

\begin{theo}\label{th:degree of difference polynomial}
Let $a(x),b(x)\in\bC[x]$ and $d,m_0$ be given by \eqref{eq:d} and \eqref{eq:m_0} if defined.
For any polynomial $g(x)\in\bC[x]$, let
\[
p(x)=a(x)g(q x)-b(q^{-1}x)g(x).
\]
Then we have
\[
\deg p(x)\left\{
     \begin{array}{ll}
       <d+\deg g(x), & \hbox{if $ (a(x),b(x))$ is degenerated and $\deg g(x)=m_0$}. \\
       =d+\deg g(x), & \hbox{otherwise.}
     \end{array}
   \right.
\]
\end{theo}
\pf It is easy to see that $\deg p(x)\le d+\deg g(x)$, and the equality holds when $\deg a(x)\neq \deg b(x)$.
Next we assume $\deg a(x)=\deg b(x)=d$ and $\deg g(x)=m$.
Then $\deg p(x)<d+\deg g(x)$ if and only if
\[
\lc(a(x)g(qx))=\lc(b(q^{-1}x)g(x)),
\]
that is, $q^{m}\lc a(x)=\lc b(q^{-1}x)$,
which happens only when $(a(x),b(x))$ is degenerated and $\deg g(x)=m=m_0$.
\qed

As we have seen in the equality \eqref{eq:difference equation}, no matter
whether a given $q$-hypergeometric term $t_k$ is summable or not, once multiplied with a suitable polynomial $r(q^k)$ on $q^k$, the product can be summable.
An upper bound on the degree of $r(x)$ can be derived as follows.
\begin{theo}
Let $t_k$ be a $q$-hypergeometric term satisfying \eqref{eq:q-hyper} and $d$ given by \eqref{eq:d}.
One can find a nonzero polynomial $r(x)$ with
$\deg r(x)\leq d+1$ such that $r(q^k)\cdot t_k$ is summable.
\end{theo}
\pf
Let $p_i(x)=a(x)(qx)^i-b(q^{-1}x)x^i$, $i\in\{0,1\}$.
From identity \eqref{eq:difference equation}, we know $p_i(q^k)t_k$ is summable.
Apparently $\deg p_i(x)\leq d+1$ by Theorem \ref{th:degree of difference polynomial}.
The conclusion follows immediately by noticing that $p_0(x)$ and $p_1(x)$ can not be zero at the same time.
\qed

Let $[p(x)]=p(x)+S_{a,b}$ denote the coset of a polynomial $p(x)$.
Given any $p(x)\in\bC[x]$ and a polynomial pair $(a(x),b(x))$,
the following \emph{$q$-polynomial reduction process} provides an algorithmic way to find a representative of $[p(x)]$ with bounded degree.

\emph{The $q$-polynomial reduction process}: For any $i\in\bN$, let
\begin{equation}\label{eq:p_s}
  p_i(x)=a(x)g_i(qx)-b(q^{-1}x)g_i(x),
\end{equation}
where $g_i(x)$ is a polynomial in $\bC[x]$ of degree $i$.
We first consider the case when $(a(x),b(x))$ is not degenerated.
By Theorem \ref{th:degree of difference polynomial}, we know
\[
\deg p_i(x)=d+i, \text{ for any } i\in\bN.
\]
Then for any polynomial $p(x)$ of degree $m$ with $m\geq d$,
it can be written by the divison algorithm as
\begin{equation}\label{eq:reduction step}
 p(x)=\sum_{k=0}^{m-d}c_{k}p_{k}(x)+\tilde{p}(x),
\end{equation}
where $c_k\in\bC$ for $0\leq k\leq m-d$ and $\tilde{p}(x)$ is a polynomial of degree less than $d$.
When $(a(x),b(x))$ is degenerated, by Theorem \ref{th:degree of difference polynomial},
\[
\deg p_i(x)=d+i, \forall i\neq m_0 \text{ and }
 \deg p_{m_0}(x)<d+m_0.
\]
Then \eqref{eq:reduction step} works well except for the polynomials of degree $d+m_0$.
Thus for any polynomial $p(x)$ of degree $m$ with $m\geq d$, we can write it as
\begin{equation}\label{eq:reduction step1}
 p(x)=\sum_{\substack{0\leq k\leq m-d\\{k}\neq m_0}}c_{k}p_{k}(x)+c_{m_0}x^{d+m_0}+\tilde{p}(x),
\end{equation}
where $c_k\in\bC$ for $0\leq k\leq m-d$ and $\tilde{p}(x)$ is a polynomial with $\deg\tilde{p}(x)<d$.
Equality \eqref{eq:reduction step} (resp. \eqref{eq:reduction step1}) is called the \emph{$q$-polynomial reduction} with respect to the polynomial pair $(a(x),b(x))$ when it is not degenerated (resp. degenerated).

Note that $p_i(x)\in S_{a,b}$.
By equalities \eqref{eq:reduction step} and \eqref{eq:reduction step1}, one can see that the quotient space $\bC[x]/S_{a,b}$ is finite-dimensional.
\begin{theo}\label{th:main q-polynomial}
Let $a(x),b(x)\in\bC[x]$, $S_{a,b}$ be the $q$-difference space corresponding to $(a(x),b(x))$ and $d,m_0$ given by \eqref{eq:d} and \eqref{eq:m_0}.
We have
\[
\bC[x]/S_{a,b}=\la [x^i]:i\in r_{a,b}\ra,
\]
where
\begin{equation}\label{eq:r_{a,b}}
r_{a,b}=\left\{
     \begin{array}{ll}
       \{0,1,\ldots,d-1,d+m_0\}, & \hbox{if $ (a(x),b(x))$ is degenerated,} \\[7pt]
       \{0,1,\ldots,d-1\}, & \hbox{otherwise.}
     \end{array}
   \right.
\end{equation}
\end{theo}

To do the $q$-polynomial reduction, we first need to construct $p_i(x)$ of the form \eqref{eq:p_s}, in which $g_i(x)$ can be any polynomial of degree $i$.
When $(a(x),b(x))$ is not degenerated, the following lemma ensures that the obtained reduced polynomial $\tilde{p}(x)$ is free of the choice of $g_i(x)$.

\begin{lem}\label{lem:unique decom not generated}
Let $(a(x),b(x))$ be a polynomial pair which is not degenerated,
and $d$ be defined by \eqref{eq:d}.
For a given $p(x)\in\bC[x]$, the $q$-polynomial reduction with different choices of $g_i(x)$ in \eqref{eq:p_s} leads to the same polynomial $\tilde{p}(x)=\sum_{i=0}^{d-1}c_ix^i$ with $c_i\in\bC$, such that $p(x)-\tilde{p}(x)\in S_{a,b}$.
\end{lem}
\proof
The existence of $\tilde{p}(x)$ is guaranteed by the $q$-polynomial reduction.
We only need to prove the uniqueness.
Suppose two polynomials $\tilde{p}(x)$ and $\hat{p}(x)$ were obtained due to different choices of $g_i(x)$.
That is,
\[
  \tilde{p}(x)+S_{a,b}=\hat{p}(x)+S_{a,b}
\]
with $\deg \tilde{p}(x)<d$ and $\deg \hat{p}(x)<d$.
So $\tilde{p}(x)-\hat{p}(x)\in S_{a,b}$.
If $\tilde{p}(x)\neq\hat{p}(x)$, then $0\leq \deg (\tilde{p}(x)-\hat{p}(x))<d$.
However Theorem~\ref{th:degree of difference polynomial} tells us that the degree of every nonzero polynomial in $S_{a,b}$ is larger than or equal to $d$, which leads to a contradiction.
This completes the proof.
\qed

\section{$q$-Rational Reduction}\label{sec:q rat red}

In this section, the $q$-polynomial reduction will be generalized to the rational case.
To this aim, we first discuss the summability on the rational multiple of a given $q$-hypergeometric term $t_k$.

Recall that by the classical $q$-Gosper algorithm \cite{Koornwinder1993, A=B}, the shift quotient of $t_k$ can be decomposed as
\begin{equation}\label{eq:q-Gosper rep}
\frac{t_{k+1}}{t_k}=\frac{a(q^k)}{b(q^k)}\cdot \frac{c(q^{k+1})}{c(q^k)},
\end{equation}
where $a,b,c$ are polynomials and
\begin{equation}\label{eq:q-Gosper cond}
\gcd(a(x),b(q^hx))=1,\quad \forall h\in\bN.
\end{equation}
We call $(a(x),b(x),c(x))$ a \emph{$q$-Gosper representation} of $t_k$ when conditions \eqref{eq:q-Gosper rep} and \eqref{eq:q-Gosper cond} hold.
Moreover, $t_k$ is summable if and only if
\begin{equation}\label{eq:q-Gosper equ}
a(x)g(qx)-b(q^{-1}x)g(x)=c(x)
\end{equation}
holds for some Laurent polynomial $g(x)\in \bC[x,x^{-1}]$.
Equation \eqref{eq:q-Gosper equ} is called the \emph{$q$-Gosper equation}.

The following theorem shows when a rational multiple $\frac{A(q^k)}{B(q^k)}t_k$ of a given $q$-hypergeometric term $t_k$ can be summable.
\begin{theo}\label{th:q(x) structure}
Let $t_k$ be a $q$-hypergeometric term and $(a(x),b(x),c(x))$ its $q$-Gosper representation.
Suppose $A(x),B(x)$ are polynomials such that $\frac{A(q^k)}{B(q^k)}t_k$ is summable.
If for any $h\in\bN$, we have
\begin{equation}\label{eq:B}
\gcd (B(x),B(q^{h+1}x))=1
\end{equation}
and
\begin{equation}\label{eq:gcd}
\gcd (B(x),a(q^{-1-h}x))=\gcd (B(x),b(xq^h))=\gcd(B(x),c(x))=1,
\end{equation}
then $B(x)\mid A(x)$.
\end{theo}
\proof
Let $\tilde{t}_k=\frac{A(q^k)}{B(q^k)}t_k$.
Then we have
\[
\frac{\tilde{t}_{k+1}}{\tilde{t}_k}=
\frac{a(q^k)B(q^k)}{b(q^k)B(q^{k+1})}
\cdot
\frac{c(q^{k+1})A(q^{k+1})}{c(q^k)A(q^k)}.
\]
Utilizing \eqref{eq:B} and \eqref{eq:gcd}, one can check that
$(a(x)B(x),b(x)B(qx),c(x)A(x))$ is a $q$-Gosper representation of $\tilde{t}_k$.
As $\tilde{t}_k$ is summable, we know the corresponding $q$-Gosper equation
\begin{equation}\label{eq:Gosper AB}
B(x)\left(a(x)g(qx)-b(q^{-1}x)g(x)\right)=c(x)A(x)
\end{equation}
holds for some Laurent polynomial $g(x)$.
By condition \eqref{eq:B}, it is easy to check that
$\gcd (B(x),x^{\ell})=1, \forall \ell\in\bN$.
Then equality \eqref{eq:Gosper AB} and condition \eqref{eq:gcd} lead to
$B(x)\mid A(x)$.
\qed

By Theorem \ref{th:q(x) structure}, if $\frac{A(q^k)}{B(q^k)}t_k$ is summable with $B(x)\nmid A(x)$ and condition \eqref{eq:B} holds, then $B(x)$ must contain some factor from $c(x)$, $a(q^{-1-h}x)$ or $b(q^hx)$ for some nonnegative integer $h$.
Next, we will show that condition \eqref{eq:B} can be further abandoned when $B(x)$ only contains factors from $c(x)$, $a(q^{-1-h}x)$ or $b(q^hx)$.

To proceed, some basic notations are needed.
Let $f(x)\in\bC[x]$ be a polynomial and $n\in\bZ$.
The \emph{shift product operator of order $n$}, denoted by $SP_n$, is defined by
\begin{equation}\label{eq:SP}
SP_{n}(f(x))=\left\{
     \begin{array}{lll}
       f(x)f(qx)\cdots f(q^{n-1}x), & \hbox{when $n>0$,} \\[7pt]
       1, & \hbox{when $n=0$,} \\[7pt]
       f(q^{n}x)f(q^{n+1}x)\cdots f(q^{-1}x), & \hbox{when $n<0$.}
     \end{array}
   \right.
\end{equation}
Suppose $(a(x),b(x))$ is a polynomial pair, and $a_1(x)$ (resp. $b_1(x)$) is any factor of $a(x)$ (resp. $b(x)$).
Then
\begin{equation}\label{eq:factor}
a(x)=\tilde{a}(x)a_1(x)\text{ and }
   b(x)=\tilde{b}(x)b_1(x),
\end{equation}
for some polynomials $\tilde{a}(x),\tilde{b}(x)\in\bC[x]$.
Let $n_1,n_2\in\bN$, denote
\begin{equation}\label{eq:shift pair}
A(x)=\tilde{a}(x)a_1(q^{-n_1}x)\text{ and }
B(x)=\tilde{b}(x)b_1(q^{n_2}x).
\end{equation}
$(A(x),B(x))$ satisfying \eqref{eq:shift pair} is called the
\emph{$(a_1(x),b_1(x))$ shift pair of order $(n_1,n_2)$ for $(a(x),b(x))$}.
The definition immediately shows that
$
\deg A(x)=\deg a(x), \deg B(x)=\deg b(x)
$
and
$
\frac{\lc B(x)}{\lc A(x)}=
\frac{\lc b(x)}{\lc a(x)}\cdot q^{n_1\cdot\deg a_1(x)+n_2\cdot\deg b_1(x)},
$
which lead to the following result.
\begin{lem}\label{lem:shift pair}
Suppose the polynomial pair $(a(x),b(x))$ satisfies condition \eqref{eq:factor} and
$(A(x),B(x))$ is the $(a_1(x),b_1(x))$ shift pair of order $(n_1,n_2)$ for $(a(x),b(x))$.
Then $(A(x),B(x))$ is degenerated if and only if
$\deg a(x)=\deg b(x)=d$,
$\frac{\lc b(x)}{\lc a(x)}=q^m$ for some $m\in\bZ$ and
\[
m+n_1\cdot\deg a_1(x)+n_2\cdot\deg b_1(x)\geq d.
\]
\end{lem}
Now we are ready to present the main theorem in this section.
\begin{theo}\label{th:q-rational decomposition}
Let $t_k$ be a $q$-hypergeometric term and $\frac{t_{k+1}}{t_k}=\frac{a(q^{ k})}{b(q^{ k})}$ for some polynomials $a(x),b(x)$.
Assume $(a(x),b(x))$ satisfies condition \eqref{eq:factor}.
Then for any nonzero $p(x)\in\bC[x]$ and nonnegative integers $n_1,n_2$, we can find another polynomial $\tilde{p}(x)$ and a $q$-hypergeometric term $T_k$ such that
\begin{equation}\label{eq:q-rational decom}
p(q^k)t_k=\frac{\tilde{p}(q^k)}{SP_{-n_1}(a_1(q^k))\cdot SP_{n_2}(b_1(q^k))}t_k
+\Delta_k(T_k),
\end{equation}
where $\tilde{p}(x)=\sum\limits_{i\in r_{A,B}}c_ix^i$, $c_i\in\bC$, $(A(x),B(x))$ is the $(a_1(x),b_1(x))$ shift pair of order $(n_1,n_2)$ for $(a(x),b(x))$, and $r_{A,B}$ is defined similarly as in \eqref{eq:r_{a,b}}.
\end{theo}
\proof
Let $s_k=\frac{t_k}{SP_{-n_1}(a_1(q^k))\cdot SP_{n_2}(b_1(q^k))}$.
It is straightforward to check that
\[
\frac{s_{k+1}}{s_k}=\frac{A(q^k)}{B(q^k)}.
\]
Let $S_{A,B}$ be the $q$-difference space corresponding to $(A(x),B(x))$.
Notice that the shift product operator $SP_n$ sends a polynomial of degree $\ell$ to a polynomial of degree $|n|\cdot\ell$ for any integer $n$.
Theorem \ref{th:main q-polynomial} guarantees that we can find a polynomial $\tilde{p}(x)=\sum\limits_{i\in r_{A,B}}c_ix^i$ with $c_i\in\bC$, such that
\begin{equation}\label{eq:main q poly eqn}
p(x)\cdot SP_{-n_1}(a_1(x))\cdot SP_{n_2}(b_1(x))-\tilde{p}(x)\in S_{A,B}.
\end{equation}
Notice that for any $r(x)=A(x)g(qx)-B(q^{-1}x)g(x)\in S_{A,B}$, we have $r(q^k)s_k=\Delta_k(B(q^{k-1})g(q^{k})s_k)$
and thus
\begin{align*}
p(q^k)t_k
& =p(q^k)\cdot SP_{-n_1}(a_1(q^k))\cdot SP_{n_2}(b_1(q^k))\cdot s_k\\
&=\frac{\tilde{p}(q^k)}{SP_{-n_1}(a_1(q^k))\cdot SP_{n_2}(b_1(q^k))}t_k+\Delta_k(T_k),
\end{align*}
for some $q$-hypergeometric term $T_k$.
\qed

From the proof of Theorem \ref{th:q-rational decomposition}, one can see that the decomposition of $p(q^k)t_k$ into the form of \eqref{eq:q-rational decom} is totally algorithmic, which will be called the \emph{$q$-rational reduction} with respect to $t_k$.

\begin{rem}\label{rem:not degenerated}
Notice that
\[
d=\max \{\deg a(x),\deg b(x)\}=\max \{\deg A(x),\deg B(x)\}.
\]
Assume that both $(a(x),b(x))$ and $(A(x),B(x))$ are not degenerated in Theorem \ref{th:q-rational decomposition}.
By Theorem \ref{th:main q-polynomial}, for any polynomial $p(x)$ with
\[
\deg p(x)+n_1\cdot\deg a_1(x)+n_2\cdot\deg b_1(x)\geq d,
\]
one can find a polynomial $\tilde{p}(x)$ with $\deg \tilde{p}(x)<d$ such that \eqref{eq:q-rational decom} holds.
\end{rem}

Identity \eqref{eq:q-rational decom} can be rewritten as
\begin{equation}\label{eq: rational summable}
\left(p(q^k)-\frac{\tilde{p}(q^k)}{SP_{-n_1}(a_1(q^k))\cdot SP_{n_2}(b_1(q^k))}\right)\cdot t_k
=\Delta_k(T_k).
\end{equation}
Here $\tilde{p}(x)$ is a polynomial with bounded degree and $n_1$,$n_2$ can be any nonnegative integer.
Thus the left hand side of the equality \eqref{eq: rational summable} is a $q$-hypergeometric term of the form $r(q^k)t_k$ with $r(x)$ being rational.
Apparently $r(q^k)t_k$ is summable but the denominator of $r(x)$ may not satisfy condition \eqref{eq:B} when $n_1>1$ or $n_2>1$.

\section{Applications}\label{sec: application}
In this section, we first show how Theorem \ref{th:q-rational decomposition} provides a uniform way to construct new $q$-identities from known ones.
When the shift quotient of $t_k$ is a rational funtion in $q^{\ell k}$ for some integer $\ell\geq 2$, we can modify the $q$-rational reduction process to obtain $\tilde{p}(x)$ of lower degree.
This enables us to obtain $q$-analogues of some Ramanujan-type series from \cite{Sun2021,Sun2020} in Subsection \ref{subsec: 2}.
At the end of this section, we take an example to illustrate that our $q$-rational reduction may lead to a rational part with smaller numerator than that obtained directly from the $q$-Gosper algorithm.

\subsection{Generating new $q$-identities}

In 1993, using the WZ method, Zeilberger \cite{Zeil1993} discovered that
\[
\sum_{k=1}^{\infty} \frac{21k-8}{k^3\binom{2k}{k}^3}=\frac{\pi^2}{6}.
\]
Following Zeilberger's work, Guillera \cite{Guillera2008} later obtained the Zeilberger-type series
\begin{equation}\label{eq:zeil 1}
\sum_{k=1}^{\infty} \frac{(3k-1)16^k}{k^3\binom{2k}{k}^3}=\frac{\pi^2}{2}.
\end{equation}
Recently, a $q$-analogue of identity \eqref{eq:zeil 1} was presented by Hou, Krattenthaler and Sun in \cite{HouKraSun2019} as
\begin{equation}\label{eq:q-zeil 1}
\sum_{k=0}^{\infty}q^{k(k+1)/2}[3k+2]_q
\cdot \frac{(q;q)_k^3(-q;q)_k}{(q^3;q^2)_k^3} =(1-q)^2\frac{(q^2;q^2)_{\infty}^4}{(q;q^2)_{\infty}^4}.
\end{equation}
Here $[n]_q$ is the $q$-analogue of $n$ given by
\[
[n]_q=\frac{1-q^n}{1-q}=\sum_{k=0}^{n-1} q^k.
\]
We also adopt the standard notation from \cite{GR2004}:
\[
(a;q)_{\infty}=\prod_{i=0}^{\infty}(1-aq^i) \text{ and }
(a;q)_n=\frac{(a;q)_{\infty}}{(aq^n;q)_{\infty}}
\]
for a nonnegative integer $n$.
We will take \eqref{eq:q-zeil 1} as an example to demonstrate how to generate new identities using the $q$-rational reduction.

\begin{prop}
The following identity is true:
\begin{equation}\label{eq: new qzeil}
\sum_{k=0}^{\infty}q^{\frac{k^2+5k+2}{2}}\frac{[3k+4]_q}
{[2k+3]_q^2}
\cdot \frac{(q;q)_k^3(-q;q)_k}{(q^3;q^2)_k^3}  =(1-q)^2\frac{(q^2;q^2)_{\infty}^4}{(q;q^2)_{\infty}^4}-(1+q).
\end{equation}
\end{prop}
\proof
For the first step, write the summand of the left hand side of \eqref{eq:q-zeil 1} as $p(q^k)\cdot t_k$, where
$p(x)=\frac{1-q^2x^3}{1-q}$ and $t_k=q^{k(k+1)/2}\cdot \frac{(q;q)_k^3(-q;q)_k}{(q^3;q^2)_k^3}$.
It is easy to check that
\[
\frac{t_{k+1}}{t_k}=\frac{a(q^{k})}{b(q^{k})}
\text{ with }
a(x)=qx(1-qx)^3(1+qx) \text{ and }b(x)=(1-q^3x^2)^3.
\]
Notice that $t_k$ is not degenerated since $\deg a(x)\neq
\deg b(x)$.
Take
\[
a_1(x)=1\text{ and }b_1(x)=\left(\frac{1-q^3x^2}{1-q}\right)^2,
\]
which are factors of $a(x)$ and $b(x)$ respectively,
and $n_1=0,n_2=1$.
Then
\[
SP_{-n_1}(a_1(q^{k}))=1
\text{ and }
SP_{n_2}(b_1(q^{k}))=[2k+3]_q^2.
\]
Let
\[
s_k=\frac{t_k}{SP_{-n_1}(a_1(q^{k}))SP_{n_2}(b_1(q^{k}))}
=\frac{t_k}{[2k+3]_q^2}.
\]
One can check that
\[\frac{s_{k+1}}{s_k}=\frac{A(q^{k})}{B(q^{k})}\]
where
$(A(x),B(x))$ is the $(a_1(x),b_1(x))$ shift pair of order $(0,1)$ for $(a(x),b(x))$.
Note that $s_k$ is not degenerated since $\deg A(x)\neq \deg B(x)$.
Then we know by Theorem \ref{th:main q-polynomial} that
\[
\bC[x]/S_{A,B}=\la [x^0],[x^1],\ldots,[x^{d-1}]\ra,
\]
where $S_{A,B}$ is defined by \eqref{eq:Sab} and
\[d=\max \{\deg A(x),\deg B(x)\}=6.\]
Notice that
$
[3k+2]_q\cdot t_k=f(q^{k}) \cdot s_k
$
with
\[
f(x)=\frac{(1-q^2x^3)(1-q^3x^2)^2}{(1-q)^3}.
\]
Let
\begin{equation}\label{eq:A,B}
p_i(x)=A(x)(q x)^i - B(q^{-1}x) x^i,\quad i\in\bN
\end{equation}
which are polynomials in $S_{A,B}$ (by taking $g(x)=x^i$).
Then
\[
p_i(q^k)s_k=\Delta_k(B(q^{k-1}) q^{ik}s_k).
\]
As $s_k$ is not degenerated, Theorem \ref{th:degree of difference polynomial} shows $\deg(p_i(x))=i+6$.
Since $\deg f(x)=7>d-1$, by the $q$-rational reduction, $f(x)$ can be rewritten as
\[
f(x)=\tilde{p}(x)-\frac{ p_0(x)}{(1-q)^3}
    -\frac{q\cdot p_1(x)}{(1-q)^3},
\]
where $\tilde{p}(x)=\frac{qx^2(1-q^4x^3)}{1-q}$, $p_0(x)$ and $p_1(x)$ are defined by \eqref{eq:A,B}.
Thus
\begin{equation}\label{eq: new qzeil1}
 p(q^k)\cdot t_k=f(q^{k}) \cdot s_k
=q^{2k+1}[3k+4]_q s_k+\Delta_k(T_k),
\end{equation}
with
\[
T_k=-q^{k(k+1)/2}(1+ q^{k+1})[2k+1]_q\cdot \frac{(q;q)_k^3(-q;q)_k}{(q^3;q^2)_k^3}.
\]
Summing identity \eqref{eq: new qzeil1} over $k$ from $0$ to $\infty$, we obtain equality \eqref{eq: new qzeil} with the help of \eqref{eq:q-zeil 1}.
\qed

Let $q\rightarrow 1$ on both sides of equality \eqref{eq: new qzeil} and then make the substitution $n=k+1$ will arrive at
\[
\sum_{n=1}^{\infty}\frac{3n+1}{(2n+1)^2}\cdot
\frac{16^n}{n^3\binom{2n}{n}^3}
=\frac{\pi^2}{2}-4,
\]
which can be seen as a rational Zeilberger-type series compared with  \eqref{eq:zeil 1}.

Equality \eqref{eq: new qzeil} is obtained by taking
\[
(a_1(x),b_1(x))=\left(1,\left(\frac{1-q^3x^2}{1-q}\right)^2\right)\text{ and }
(n_1,n_2)=(0,1).
\]
If we take
\[
(a_1(x),b_1(x))=\left(1,\frac{1-q^3x^2}{1-q}\right)
\text{ and }
(n_1,n_2)=(0,2),
\]
similar discussions will arrive at
\[
\sum_{k=0}^{\infty}\frac{P(q^k)\cdot q^{k(k+1)/2}}
{[2k+3]_q[2k+5]_q}
\cdot \frac{(q;q)_k^3(-q;q)_k}{(q^3;q^2)_k^3}  =(1-q)^2\frac{(q^2;q^2)_{\infty}^4}{(q;q^2)_{\infty}^4}
-\frac{q^2(1+q)}{1+q+q^2},
\]
where
$P(x)=\frac{1}{(1-q)^2}(1+q-(-q^3 + q^4 + q^5 + q^6) x^2-(q^2 + q^3 - q^4 - q^5 + q^6 + q^7) x^3-(q^4 + q^5 - q^6 - q^7) x^4+2 q^8 x^5)$,
which is a $q$-analogue of
\begin{equation}\label{eq:zeil shift}
\sum_{n=1}^{\infty}\frac{4 n^2+ 7 n-3 }{(2n+1)(2n+3)}\cdot
\frac{16^n}{n^3\binom{2n}{n}^3}
=\frac{\pi^2}{6}-\frac{4}{9}.
\end{equation}
Remark \ref{rem:not degenerated} ensures that there are infinitely many other choices of $(a_1(x),b_1(x))$ and $(n_1,n_2)$, in principle any choice will lead to a $q$-series that is a $q$-analogue of some rational Zeilberger-type series.

It should be mentioned that, with the help of the extended Zeilberger's algorithm, Hou and Li \cite{HouLi2021} also provided a systematic method to construct new hypergeometric series.
But they only focused on those with the summand of the form $\frac{a(k)}{b(k)}\cdot t_k$, where $t_k$ is hypergeometric, $a(k)$, $b(k)$ are polynomials with $\gcd(b(k),b(k+1+h))=1$ for any $h\in\bN$.
Thus \eqref{eq:zeil shift} and infinitely many other similar identities will be missed.
\begin{rem}
In general, $p_i(x)\in S_{A,B}$ is of the form
\[
p_i(x)=A(x)g(qx)-B(q^{-1}x)g(x),
\]
where $g(x)\in\bC[x]$ can be any polynomial of degree $i$.
Note that when $(A(x),B(x))$ is not degenerated, Lemma \ref{lem:unique decom not generated} ensures the generated identities are free of the choice of $g(x)$.
In this paper, we always take $g(x)=x^i$.
\end{rem}
\subsection{When a shift quotient is a polynomial in $q^{\ell k}$}\label{subsec: 2}
As we have seen, once a summation identity with the summand of the  form $p(q^k)t_k$ is given, the $q$-rational reduction may lead to new equalities.
Here $p(x)\in\bC[x]$ and $t_k$ is a $q$-hypergeometric term.
Sometimes, the shift quotient of $t_k$ is a rational function not only in $q^k$ but also in $q^{\ell k}$ for some $\ell\in\bN$ with $\ell\geq 2$.
In this case, we will show how to modify the $q$-rational reduction to decrease the dimension of the corresponding quotient space.

We will take a series for $\pi^{-1}$ as an example, the history of which was begun by G. Bauer \cite{Bauer1859} in 1859 with the discovery of the following identity
\begin{equation*}\label{eq:classical R -64}
 \sum_{k=0}^{\infty}(4k+1)\frac{\binom{2k}{k}^3}{(-64)^k}
 =\frac{2}{\pi}.
\end{equation*}

In Bauer's series, the summand can be seen as a product of the polynomial $(4k+1)$ and a hypergeometric term $\frac{\binom{2k}{k}^3}{(-64)^k}$.
Recently the following four rational-type Bauer's series, with the polynomial part replaced by a rational function, were obtained by Z.-W. Sun in \cite{Sun2020} utilizing Gosper's algorithm.
\begin{align}
&\sum_{k=0}^{\infty}
\frac{k(4 k-1)}{(2k-1)^2}\frac{\binom{2k}{k}^3}{(-64)^k}
=-\frac{1}{\pi},\label{eq:rat Bauer (2k-1)^2}\\
&\sum_{k=0}^{\infty}
\frac{( 4 k-1)}{(2k-1)^3}\frac{\binom{2k}{k}^3}{(-64)^k} =\frac{2}{\pi},
\label{eq:rat Bauer (2k-1)^3}\\
&\sum_{k=0}^{\infty}\frac{(4k+1)
 \binom{2k}{k}^3}{(2k-1)(k+1)(-64)^k}
 =-\frac{4}{\pi},\label{eq:rat Bauer (2k-1)(k+1)}\\
&\sum_{k=0}^{\infty}
 \frac{(2k+1)(4k+3)\binom{2k}{k}^3}{(k+1)^2(-64)^k} =\frac{8}{\pi} \label{eq:rat Bauer (k+1)^2}.
\end{align}

V.J.W. Guo \cite{Guo2020} pointed out that a $q$-analogue of Bauer's formula
can be written as
\begin{equation}\label{eq:q-Bauer}
\sum_{k=0}^{\infty}(-1)^kq^{k^2}[4k+1]_q
\frac{(q;q^2)_k^3}{(q^2;q^2)_k^3}=
\frac{(q;q^2)_{\infty}(q^3;q^2)_{\infty}}{(q^2;q^2)_{\infty}^2}.
\end{equation}

By the $q$-Gosper algorithm and identity \eqref{eq:q-Bauer}, Q.-H. Hou and Z.-W. Sun \cite{HouSun2021} obtained $q$-analogues of \eqref{eq:rat Bauer (2k-1)^2} and \eqref{eq:rat Bauer (2k-1)^3} as follows:
\begin{equation} \label{eq:(2k-1)^2}
\sum_{k=0}^{\infty}(-1)^kq^{k^2}
\frac{[2k]_q([4k]_q-1)}
{([2k]_q-1)^2}
\cdot \frac{(q;q^2)_k^3}{(q^2;q^2)_k^3} =-\frac{(q;q^2)_{\infty}(q^3;q^2)_{\infty}}
{(q^2;q^2)_{\infty}^2}.
\end{equation}
\begin{equation}\label{eq:(2k-1)^3}
\sum_{k=0}^{\infty}(-1)^kq^{k^2+2k}
\frac{([4k]_q-1)}
{([2k]_q-1)^3}
\cdot \frac{(q;q^2)_k^3}{(q^2;q^2)_k^3} =\frac{(q;q^2)_{\infty}(q^3;q^2)_{\infty}}
{(q^2;q^2)_{\infty}^2}.
\end{equation}

Let $A_k(q)$ be the summand on the left hand side of \eqref{eq:q-Bauer}.
The critical step in the proof of the above two equalities is the construction of another $q$-hypergeometric term $B_k(q)$ such that $A_k(q)+B_k(q)$ is summable.
However, how to find the $B_k(q)$ by the $q$-Gosper algorithm is empirical.

In Section \ref{sec:q rat red}, we introduce the $q$-rational reduction which offers a mechanical pattern to find the $B_k(q)$.
Next we will show that $q$-analogues of the four rational-type Bauer's identities can be obtained uniformly by the $q$-rational reduction.
More precisely, we can not only reprove \eqref{eq:(2k-1)^2} and \eqref{eq:(2k-1)^3}, but also provide $q$-analogues of \eqref{eq:rat Bauer (2k-1)(k+1)} and \eqref{eq:rat Bauer (k+1)^2} as follows.

For simplicity, we define
\begin{equation*}
SP_{n}^{(2)}(f(x))=\left\{
     \begin{array}{lll}
       f(x)f(q^2x)\cdots f(q^{2n-2}x), & \hbox{when $n>0$,} \\[7pt]
       1, & \hbox{when $n=0$,} \\[7pt]
       f(q^{2n}x)f(q^{2n+2}x)\cdots f(q^{-2}x), & \hbox{when $n<0$.}
     \end{array}
   \right.
\end{equation*}
\begin{prop}\label{eq:all  q-bauer}
The following identities are true:
\begin{align}
& \sum_{k=0}^{\infty}(-1)^kq^{k^2+2k}\frac{[4k+1]_q}
{([2k]_q-1)[2k+2]_q}
\cdot \frac{(q;q^2)_k^3}{(q^2;q^2)_k^3} =-\frac{(q;q^2)_{\infty}(q^3;q^2)_{\infty}}
{(q^2;q^2)_{\infty}^2}. \label{eq:(2k-1)(k+1)} \\
&\sum_{k=0}^{\infty}(-1)^kq^{k^2}\frac{P(q^{2k})}
{(1-q)^2[2k+2]_q^2}
\cdot \frac{(q;q^2)_k^3}{(q^2;q^2)_k^3} =\frac{(q;q^2)_{\infty}(q^3;q^2)_{\infty}}
{(q^2;q^2)_{\infty}^2},\label{eq:(2k+1)^2}
\end{align}
where $P(x)=
1+(q -2 q^{2})x - (q  + q^2 - q^4 )x^2+ (2 q^3  - q^4 )x^3$.
\end{prop}
\proof
Identities \eqref{eq:(2k-1)^2}--\eqref{eq:(2k+1)^2} can be proved similarly.
To compare with the $q$-Gosper algorithm, we will prove identity \eqref{eq:(2k-1)(k+1)} in detail.

The first step is to write the summand of the left hand side of \eqref{eq:q-Bauer} as $[4k+1]_q\cdot t_k$, where $t_k=(-1)^kq^{k^2}\cdot \frac{(q;q^2)_k^3}{(q^2;q^2)_k^3}$ is a $q$-hypergeometric term and
$[4k+1]_q=\frac{1-q^{4k+1}}{1-q}$ is a polynomial in $q^{2k}$.
It is easy to check that
\[
\frac{t_{k+1}}{t_k}=\frac{a(q^{2k})}{b(q^{2k})}
\text{ with }
a(x)=(-q)x(1-qx)^3 \text{ and }b(x)=(1-q^2x)^3.
\]
Motivated by the structure of the summand on the left hand side of \eqref{eq:rat Bauer (2k-1)(k+1)}, we will take
\[
a_1(x)=\frac{q(1-qx)}{1-q}\text{ and }b_1(x)=\frac{1-q^2x}{1-q},
\]
which are factors of $a(x)$ and $b(x)$ respectively.
Let
$\tilde{a}(x)=a(x)/a_1(x)$ and $\tilde{b}(x)=b(x)/b_1(x)$.
Notice that
\[
SP_{-1}^{(2)}(a_1(q^{2k}))=[2k]_q-1
\text{ and }
SP_{1}^{(2)}(b_1(q^{2k}))=[2k+2]_q.
\]
Let
\[
s_k=\frac{t_k}{([2k]_q-1)[2k+2]_q}.
\]
One can check that
\[\frac{s_{k+1}}{s_k}=\frac{A(q^{2k})}{B(q^{2k})}\]
where
$
A(x)=\tilde{a}(x)a_1(q^{-2}x)\text{ and }
B(x)=\tilde{b}(x)b_1(q^{2}x).
$
For any $g(x)\in\bC[x]$, note that
\[
\left(A(q^{2k})g(q^{2k+2})-B(q^{2k-2})g(q^{2k})\right)s_k
=\Delta_k(B(q^{2k-2})g(q^{2k})s_k),
\]
Similar to \eqref{eq:Sab}, we denote
\[
S_{A,B}^{(2)}:=\{A(x)g(q^2x)-B(q^{-2}x)g(x)\mid g(x)\in\bC[x]\}.
\]
Then $g(q^{2k}) s_k$ is summable for any polynomial $g(x)\in S_{A,B}^{(2)}$.
Note that $s_k$ is not degenerated since $\deg A(x)\neq \deg B(x)$.
Let
\begin{equation}\label{eq:A,B,2}
p_i(x)=A(x)(q^2 x)^i - B(q^{-2}x) x^i,\quad i\in\bN,
\end{equation}
which are polynomials in $S_{A,B}^{(2)}$ (by taking $g(x)=x^i$).
Similar discussions as in the proof of Theorem \ref{th:degree of difference polynomial} and Theorem \ref{th:main q-polynomial} lead to $\deg(p_i(x))=i+d$ and
\[
\bC[x]/S_{A,B}^{(2)}=\la [x^0],[x^1],\cdots,[x^{d-1}]\ra,
\]
where  $d=\max\{\deg A(x),\deg B(x)\}=4$.
Notice that
$
[4k+1]_q\cdot t_k=f(q^{2k}) \cdot s_k
$
 with
$
f(x)=\frac{1}{(1-q)^3}\cdot (1-qx^2)(q-x)(1-q^2x).
$
Since $\deg f(x)=4>d-1$, $f(q^{2k})$ can be rewritten as
\[
f(q^{2k})=-\frac{q\cdot p_0(q^{2k})}{(1-q)^3}-q^{2k}[4k+1]_q,
\]
by the $q$-polynomial reduction.
Thus
\begin{equation}\label{eq:cal (2k-1)^2}
 [4k+1]_q\cdot t_k=f(q^{2k}) \cdot s_k
=-q^{2k}[4k+1]_q s_k+\Delta_k(T_k),
\end{equation}
with
\[
T_k=(-1)^{k+1} q^{k^2+1}\frac{[2k]_q^2[2k+2]_q}{([2k]_q-1)[2k+2]_q}\cdot
\frac{(q;q^2)_k^3}{(q^2;q^2)_k^3}.
\]
Summing identity \eqref{eq:cal (2k-1)^2} over $k$ from $0$ to $\infty$, we will obtain equality \eqref{eq:(2k-1)(k+1)}
by \eqref{eq:q-Bauer}.
\qed

\subsection{Rational part with the numerator of lower degree}
As we have seen, compared to the $q$-Gosper algorithm, our $q$-rational reduction process can not only determine the structure of the denominator of the rational part of the series but also characterise the degree of the numerator.
Next, we will show our $q$-rational reduction process may lead to rational part with the numerator of lower degree than those obtained by the $q$-Gosper algorithm.

Utilizing the $q$-WZ method, Guo and Liu \cite{GuoLiu2018} proved
\begin{equation}\label{eq:q-rama 256}
 \sum_{k=0}^{\infty}
 q^{k^2} [6k+1]_q
 \frac{(q;q^2)_k^2(q^2;q^4)_k}{(q^4;q^4)_k^3}=
 (1+q)\frac{(q^2;q^4)_{\infty}(q^6;q^4)_{\infty}}
 {(q^4;q^4)_{\infty}^2}
\end{equation}
as a $q$-analogue of Ramanujan's series \cite{Ramanujan1914}
\begin{equation*}\label{eq:rama 256}
\sum_{k=0}^{\infty} (6k+1)\frac{\binom{2k}{k}^3}{256^k}
=\frac{4}{\pi}.
\end{equation*}
In \cite{Sun2020}, Z.-W. Sun presented the following identities:
 \begin{equation}\label{eq: 256 (2k-1)^2}
\sum_{k=0}^{\infty} \frac{(12k^2-1)\binom{2k}{k}^3}{(2k-1)^2 256^k}
=-\frac{2}{\pi}
\end{equation}
and
\begin{equation}\label{eq:256 (2k-1)^3}
\sum_{k=0}^{\infty}
\frac{k( 6 k-1)}{(2k-1)^3}\frac{\binom{2k}{k}^3}{(256)^k} =\frac{1}{2\pi}.
\end{equation}

Recently Q.-H. Hou and Z.-W. Sun \cite[Theorem 1.2]{HouSun2021} proved
\begin{equation}\label{eq:HS q-256 (2k-1)^2}
\sum_{k=0}^{\infty} \frac{P\cdot q^{k^2}}{(1-q)^3([2k]_q-1)^2}
\cdot
\frac{(q;q^2)_k^2(q^2;q^4)_k}{(q^4;q^4)_k^3}
=2q(1+q)\frac{(q^2;q^4)_{\infty}(q^6;q^4)_{\infty}}
 {(q^4;q^4)_{\infty}^2},
\end{equation}
where $P$ denotes
\[
q^{12k+1}-3q^{10k+2}+3(2q^2-1)q^{8k+1}-(3q^4-1)q^{6k}
+3q^{4k+1}-3q^{2k+2}+2q^3-q.
\]
Identity \eqref{eq:HS q-256 (2k-1)^2} can be seen as a $q$-analogue of equality \eqref{eq: 256 (2k-1)^2}.
However, no $q$-analogue of \eqref{eq:256 (2k-1)^3} was given.

Next, we will provide another $q$-analogue of \eqref{eq: 256 (2k-1)^2} with a ``smaller" $P$ and a $q$-analogue of \eqref{eq:256 (2k-1)^3} by the $q$-rational reduction.
\begin{prop}
The following identities are $q$-analogues of identities \eqref{eq: 256 (2k-1)^2} and \eqref{eq:256 (2k-1)^3} respectively:
\begin{equation}\label{eq:(2k-1)(4k-2)}
\sum_{k=0}^{\infty}\frac{P_1\cdot q^{k^2-4}}{(1 - q)^2[2k-1]_q[4k-2]_q}
\cdot \frac{(q;q^2)_k^2(q^2;q^4)_k}{(q^4;q^4)_k^3}
=\frac{(q^2;q^4)_{\infty}(q^6;q^4)_{\infty}}
 {(q^4;q^4)_{\infty}^2}
\end{equation}
with
$
P_1=-q^2 +  q^{4k+2} + q^{6 k+3}+q^{8k} - q^{8 k+2} -q^{10 k+1},
$
and
\begin{equation}\label{eq:(2k-1)^3 256^k}
\sum_{k=0}^{\infty}\frac{P_2\cdot q^{k^2-7}(1+q)}{(1 - q)^2[2k-1]_q^2[4k-2]_q}
\cdot \frac{(q;q^2)_k^2(q^2;q^4)_k}{(q^4;q^4)_k^3}
=\frac{(q^2;q^4)_{\infty}(q^6;q^4)_{\infty}}
 {(q^4;q^4)_{\infty}^2}
\end{equation}
with
$
P_2=q^3  - q^{4 k+3} - q^{ 6 k+4} - 2 q^{ 8 k+1} + 2 q^{ 8 k+3}+ q^{10 k}.
$
\end{prop}
\pf
We will start with the equality \eqref{eq:q-rama 256}, the summand on the left hand side of which can be written as $[6k+1]_q\cdot t_k$ with
$
t_k=q^{k^2}\cdot\frac{(q;q^2)_k^2(q^2;q^4)_k}{(q^4;q^4)_k^3}
$
and
$
\frac{t_{k+1}}{t_k}=\frac{a(q^{2k})}{b(q^{2k})}.
$
Here
$
a(x)=qx(1-qx)^2(1-q^2x^2)\text{ and }
b(x)=(1-q^4x^2)^3.
$
Let $a_1(x)=\frac{1-qx}{1-q}\cdot \frac{1-q^2x^2}{1-q}$ and $b_1(x)=1$.
Then
$
SP_{-1}^{(2)}(a_1(q^{2k}))=[2k-1]_q[4k-2]_q
$
and
$
SP_{1}^{(2)}(b_1(q^{2k}))=1.
$
Denote
\[
s_k=\frac{q^{k^2}}{[2k-1]_q[4k-2]_q}
\cdot
    \frac{(q;q^2)_k^2(q^2;q^4)_k}{(q^4;q^4)_k^3}
\]
and
$
f(x)=\frac{(1 - q^{-1}x) (1 - q^{-2}x^2) (1 - q x^3)}{(1 - q)^3}.
$
We have
$
[6k+1]_q\cdot t_k=f(q^{2k})\cdot s_k,
$
and
$
\frac{s_{k+1}}{s_k}=\frac{A(q^{2k})}{B(q^{2k})}
$
with
$
A(x)=x (q - x)^2 (q + x) (1 - q x)\text{ and }
B(x)=q^2 (1 - q^4 x^2)^3.
$
Apparently $s_k$ is not degenerated and $d=6$.
Take
\[
p_i(x)=A(x) (q^2 x)^i-B(q^{-2}x) x^i.
\]
Theorem \ref{th:degree of difference polynomial} guarantees
$\deg p_0(x)=\deg f(x)=6$.
By the $q$-polynomial reduction, we find
\begin{equation}\label{eq:p(q^{2k})}
f(q^{2k})=-\frac{p_0(q^{2k})}{(1 - q)^3 q^4}+\frac{(1 + q) P_1}{(1 - q)^2 q^4}.
\end{equation}
Multiplied by $s_k$ on both sides of \eqref{eq:p(q^{2k})}, we have
\begin{equation}\label{eq:q^{k^2} [6k+1]_q}
  q^{k^2} [6k+1]_q
 \frac{(q;q^2)_k^2(q^2;q^4)_k}{(q^4;q^4)_k^3}=
\Delta_k(T_k)+ \frac{(1 + q) P_1}{(1 - q)^2 q^4}\cdot s_k
\end{equation}
with
\[
T_k=
\frac{-q^{1 + k^2} (1 - q^{4 k})^3}
     {(1 - q) (q - q^{2 k})^2 (q + q^{2 k})}
\cdot\frac{(q;q^2)_k^2(q^2;q^4)_k}{(q^4;q^4)_k^3}.
\]
Summing \eqref{eq:q^{k^2} [6k+1]_q} over $k$ from $0$ to $\infty$, one can obtain \eqref{eq:(2k-1)(4k-2)} by the equality \eqref{eq:q-rama 256}.
If we take $a_1(x)=\left(\frac{1-qx}{1-q}\right)^2\cdot \frac{1-q^2x^2}{1-q}$ and $b_1(x)=1$,
similar discussions will lead to \eqref{eq:(2k-1)^3 256^k}.
\qed

\noindent \textbf{Acknowledgments.}
This work was supported by the National Natural Science Foundation of China (No. 12101449, 11701419, 11871067) and the Natural Science Foundation of Tianjin, China (No. 19JCQNJC14500).

\end{document}